\documentclass{article}
\begin{document}
\centerline{\textbf{Formulas for the approximation of the complete Elliptic Integrals}}
\[
\]
\centerline{\bf Nikos Bagis\rm}
\centerline{Department of Informatics} 
\centerline{Aristotle University  Thessaloniki Greece.}
\centerline{e-mail: nikosbagis@hotmail.gr}
\[
\]
\begin{quote}
\begin{abstract}
In this article we give evaluations of the two complete elliptic integrals $K$ and $E$ in the form of Ramanujans type-$\pi$ formulas. The result is a formula for $\Gamma(1/4)^2\pi^{-3/2}$ with accuracy about 120 digits per term.     
\end{abstract}

\bf keywords: \rm{elliptic functions; singular modulus; Ramanujan; Legendre functions; evaluations; constants}

\end{quote}

\section{Elliptic singular moduli}

It is known that (see [1],[3]) 
\begin{equation}
K(x)=\int^{\pi/2}_{0}\frac{d\theta}{\sqrt{1-x^2\sin^2(\theta)}}=\frac{\pi}{2}{}_2F_1\left(\frac{1}{2},\frac{1}{2};1;x^2\right)\end{equation}
is the complete elliptic integral of the first kind. The function $k_r$ is called elliptic singular moduli and defined from the equation 
\begin{equation}
\frac{K\left(\sqrt{1-k^2_r}\right)}{K(k_r)}=\sqrt{r}
\end{equation}
Also it is known that if $r\in\bf Q^{*}_{+}\rm$, the $k_r$ are algebraic numbers.\\
For $r\in\bf N \rm$ we set $K[r]=K(k_r)$. $K[r]$ could be expressed in terms of products of $\Gamma$ functions, algebraic numbers and powers of $\pi$. In time it became obvious that the best way to express the numbers $K[r]$ most concisely was to use the function 
\begin{equation}
b(p)=\frac{\Gamma^2(p)}{\Gamma(2p)}\sqrt{\tan(p\pi)}
\end{equation}
It is also known that if $N=n^2m$, where $n$ and $m$ are positive integers then
\begin{equation}
K[n^2m]=M_n(m)K[m] , 
\end{equation}
where $M_n(m)$ is algebraic. The following formulas for some $M_n(m)$ are known.
\begin{equation}
M_2(m)=\frac{1+k'_m}{2}
\end{equation}
\begin{equation}
27M^4_3(m)-18M^2_3(m)-8(1-2k^2_m)M_3(m)-1=0
\end{equation}
\begin{equation}
(5M_5(m)-1)^5(1-M_5(m))=256k^2_m(1-k^2_m)M_5(m)
\end{equation}
These formulas for finding $K[4r]$, $K[9r]$ and $K[25r]$ depend only on knowing $k_r$.\\
Also we consider the complete elliptic integral of the second kind, which is 
\begin{equation}
E(x)=\frac{\pi}{2}{}_2F_1\left(\frac{-1}{2},\frac{1}{2};1;x^2\right)
\end{equation}
and related with $K(x)$ from the relation
\begin{equation}
E(k_r)=\frac{K(k_r)}{\sqrt{r}}\left(\frac{\pi}{3K(k_r)^2}-a(r)\right)+K(k_r).
\end{equation}
The function $a(r)$ is called elliptic alpha function (see [4]). \\
We will use the elliptic functions theory to evaluate values of $K(k_r)$ and $E(k_r)$ in high precision using Ramanujan's type-$\pi$ formulas, but now the constant will be not $\pi$ but   
\begin{equation}
\frac{1}{\pi}b\left(\frac{1}{4}\right)=\frac{\Gamma\left(\frac{1}{4}\right)^2}{\pi^{3/2}}
\end{equation}
the precision of the application formula, which is our more interesting result in this paper is an about 120 digits per term.\\  
Our methods consists Legendre functions, and we not use the function $a(r)$.  

\section{Legendre polynomials and the formula}

The Legendre $P$ function is defined by
\begin{equation}
P^{\mu}_{\nu}(z)=\frac{1}{\Gamma(1-\nu)}\left(\frac{z+1}{1-z}\right)^{\nu/2} {}_2F_{1}\left(-\mu,\mu+1;1-\nu;\frac{1-z}{2}\right)
\end{equation}
Set
$$
\phi(z)={}_2F_{1}\left(-\mu,\mu+1;1-\nu;z\right)=\left(\frac{z}{1-z}\right)^{\nu/2}\Gamma(1-\nu)P^{\mu}_{\nu}(1-2z)
$$
Then derivating $\phi$ we have
$$
\phi'(z)=\frac{1}{2(1-z)z}\left(\frac{z}{1-z}\right)^{\nu/2} \Gamma(1-\nu)\times
$$
\begin{equation}
\times[\left(-1-\mu+\nu+2(1+\mu)z\right)P^{\mu}_{\nu}(1-2 z)+(1+\mu-\nu)P^{1+\mu}_{\nu}(1-2z)]
\end{equation}
If we assume that
\begin{equation}
\sum^{\infty}_{n=0}\frac{\left(-\mu\right)_n\left(1+\mu\right)_n}{\left(1-\nu\right)_n}\frac{z^n}{n!}(\alpha n+\beta)=g
\end{equation}
then
$$\beta \phi(z)+\alpha z\phi'(z)=g$$
From (11),(12) and (13) we have\\ 
\textbf{Theorem 1.}\\
If 
\begin{equation}
\alpha=\frac{2(-1+z)}{-1-\mu+\nu+2z+2\mu z}
\end{equation}
then
\small
\begin{equation}
\sum^{\infty}_{n=0}\frac{\left(-\mu\right)_n\left(1+\mu\right)_n}{\left(1-\nu\right)_nn!}z^n(\alpha n+1)
=\frac{(-1-\mu+\nu)\left(\frac{z}{1-z}\right)^{\nu/2} \Gamma(1-\nu)P^{1+\mu}_{\nu}(1-2z)}{-1-\mu+\nu+2(\mu+1)z}
\end{equation}
\normalsize
\[
\]
It is known (see [1]), that 
\begin{equation}
P^{(-1/2)}_{0}(1-2z)={}_2F_{1}\left(\frac{1}{2},\frac{1}{2};1;z\right)
\end{equation}
hence if we set $\mu=-3/2$ and $\nu=0$, then we have\\ 
\textbf{Proposition 1.}
\begin{equation}
\sum^{\infty}_{n=0}\frac{\left(\frac{3}{2}\right)_n\left(\frac{-1}{2}\right)_n}{(n!)^2}(k_r)^{2n}\left[-4(1-k^2_r)n+1-2k^2_r\right]=\frac{2K(k_r)}{\pi}=2\vartheta^2_3(q)
\end{equation}
where $k'_{r}=\sqrt{1-k^2_r}$, $q=e^{-\pi\sqrt{r}}$. 
\[
\]
The result of the Proposition 1 is not trivial since the $\vartheta_3$-function can be evaluated from the identity
\begin{equation}  
\vartheta_3(q)=\sum^{\infty}_{n=-\infty}q^{n^2} 
\end{equation}
in which the two constants $e$ and $\pi$ involved.
\[
\]  
We know that 
\begin{equation}
k_{4r}=\frac{1-k'_r}{1+k'_{r}}
\end{equation}
Hence
$$K[16r]=\frac{1+k'_{4r}}{2}K[4r]=\frac{1+k'_{4r}}{2}\frac{1+k'_{r}}{2}K[r]$$
But
$$k'_{4r}=\sqrt{1-k^2_{4r}}=\sqrt{1-\left(\frac{1-k'_r}{1+k'_r}\right)^2}=\frac{\sqrt{(1+k'_r)^2-(1-k'_{r})^2}}{1+k'_r}=$$
or
\begin{equation}
k'_{4r}=\frac{2\sqrt{k'_r}}{1+k'_{r}}
\end{equation}
Hence 
$$K[16r]=\frac{1+k'_r+2\sqrt{k'_{r}}}{4}K[r]$$
or 
\begin{equation}
K[16r]=\left(\frac{1+\sqrt{k'_r}}{2}\right)^2K[r]
\end{equation}
Setting $r\rightarrow 4 r$ we get
$$
K[64r]=\left(\frac{1+\sqrt{k'_{4r}}}{2}\right)^2K[4r]
$$
or\\
\textbf{Lemma.}\\
If $r>0$, then
\begin{equation}
K[64r]=\frac{\left(\sqrt{1+k'_r}+\sqrt{2\sqrt{k'_r}}\right)^2}{8}K[r]
\end{equation}
\[
\]
\section{Applications}
Set
\begin{equation}
p=2+216\cdot 5^{1/4}-96\cdot 5^{3/4}
\end{equation}
then
\begin{equation}
k_{100}=\frac{2-\sqrt{p}}{2+\sqrt{p}}\textrm{ and  }k'_{100}=\frac{2\sqrt{2}p^{1/4}}{2+\sqrt{p}}
\end{equation}
From the duplication formula is
$$k_{400}=\left(\frac{\sqrt{2}-p^{1/4}}{\sqrt{2}+p^{1/4}}\right)^2\textrm{ and }k'_{400}=\frac{2^{7/3}p^{1/8}\sqrt{2+p^{1/2}}}{(\sqrt{2}+p^{1/4})^2}$$
$$
k_{1600}=\frac{\left(\sqrt{2}+p^{1/4}\right)^2-2\cdot 2^{3/4} p^{1/8}  \sqrt{2+\sqrt{p}}}{\left(\sqrt{2}+p^{1/4}\right)^2+2\cdot 2^{3/4}p^{1/8}\sqrt{2+\sqrt{p}}}
$$
\scriptsize
$
k_{6400}=w=$
$$
=\frac{2-2\cdot 2^{5/8} \left(2+\sqrt{p}\right)^{1/4} \sqrt{2 \sqrt{2}+4 p^{1/4}+\sqrt{2} \sqrt{p}} p^{1/16}+2\cdot 2^{3/4} \sqrt{2+\sqrt{p}} p^{1/8}+2 \sqrt{2} p^{1/4}+\sqrt{p}}{2+2\cdot 2^{5/8} \left(2+\sqrt{p}\right)^{1/4} \sqrt{2\cdot \sqrt{2}+4 p^{1/4}+\sqrt{2} \sqrt{p}} p^{1/16}+2\cdot 2^{3/4} \sqrt{2+\sqrt{p}} p^{1/8}+2 \sqrt{2} p^{1/4}+\sqrt{p}}
$$
\normalsize
Also from (22) we have 
$$
K[6400]=\frac{1}{8}\left(\sqrt{1+\frac{2 \sqrt{2} p^{1/4}}{2+\sqrt{p}}}+2^{7/8} \left(\frac{p^{1/4}}{2+\sqrt{p}}\right)^{1/4}\right)^2 K[100]
$$
But it is known that 
$$K[100]=\frac{4+2\sqrt{5}+\sqrt{2}(3+2\cdot 5^{1/4})}{80}b\left(\frac{1}{4}\right)$$
hence we get an about 120 digits per term formula for $\frac{1}{\pi}b(1/4)$:
\small
\[
\]
$$
\frac{1}{8}\left[4+2\sqrt{5}+\sqrt{2}(3+2\cdot 5^{1/4})\right]^{-1}\left[\sqrt{1+\frac{2 \sqrt{2} p^{1/4}}{2+\sqrt{p}}}+2^{7/8} \left(\frac{p^{1/4}}{2+\sqrt{p}}\right)^{1/4}\right]^{-2}\times
$$
$$
\times\sum^{\infty}_{n=0}\frac{\left(\frac{3}{2}\right)_n\left(\frac{-1}{2}\right)_n}{(n!)^2}(w)^{2n}\left[-2(1-w^2)n-w^2+1/2\right]=\frac{\Gamma\left(\frac{1}{4}\right)}{\pi^{3/2}}
\eqno{:(a)}
$$
\[
\]
\normalsize
\[
\]
The evaluation of $E(k_r)/\pi$ follows if we use the formula 
\begin{equation}
P_{1/2}(1-2z)=\frac{2}{\pi}[2E(z)-K(z)] , 
\end{equation}
Then one can arrive with the same method as in Proposition 1, to\\    
\textbf{Proposition 2.}
\begin{equation}
\frac{4E(k_r)}{\pi}=\frac{2K(k_r)}{\pi}+\sum^{\infty}_{n=0}\frac{\left(\frac{1}{2}\right)^2_n}{(n!)^2}(k_r)^{2n}[4(1-k^2_r)n+1-2k^2_r]
\end{equation}
\[
\]

\newpage

\centerline{\bf References}\vskip .2in

[1]: M.Abramowitz and I.A.Stegun, 'Handbook of Mathematical Functions'. Dover Publications

[2]: B.C.Berndt, 'Ramanujan`s Notebooks Part II'. Springer Verlag, New York (1989)

[3]: B.C.Berndt, 'Ramanujan`s Notebooks Part III'. Springer Verlag, New York (1991) 

[4]: J.M. Borwein and P.B. Borwein, 'Pi and the AGM'. John Wiley and Sons, Inc. New York, Chichester, Brisbane, Toronto, Singapore (1987). 

[5]: I.S. Gradshteyn and I.M. Ryzhik, 'Table of Integrals, Series and Products'. Academic Press (1980).
 
[6]: E.T.Whittaker and G.N.Watson, 'A course on Modern Analysis'. Cambridge U.P. (1927)

[7]: I.J.Zucker, 'The summation of series of hyperbolic functions'. SIAM J. Math. Ana.10.192(1979)

[8]: Bruce C. Berndt and Heng Huat Chan, 'Eisenstein Series and Approximations to Pi'. Page stored in the Web.

[9]: D. Broadhurst, 'Solutions by radicals at Singular Values $k_N$ from New Class Invariants for $N\equiv3\;\; mod\;\; 8$'. arXiv:0807.2976v3(math-phy).

\end{document}